\documentclass[12pt]{amsart}

\usepackage{amsmath}
\usepackage{upref}
\usepackage{latexsym}
\usepackage{amsfonts}
\usepackage{isolatin1}

\newtheorem{prp}{Proposition}[section]

\newtheorem{thm}{Theorem}[section]

\theoremstyle{definition}
\newtheorem{rem}{Remark}

\theoremstyle{remark}

\newcommand{\ovprt}{\overline{\partial}}
\newcommand{\ovli}{\overline}

\title{Magnetic Schr\"odinger operators and the $\ovprt $-equation}
\author{Friedrich Haslinger}

\begin{document}

\maketitle


\begin{abstract}                        
In this paper we characterize compactness of the canonical solution operator
to $\ovprt $ on weigthed $L^2$ spaces on $\mathbb C.$
For this purpose we consider certain Schr\"odinger
operators with magnetic fields and use a condition which is equivalent to the
property that these operators have compact resolvents.
We also point out what are the obstructions in the case of several complex variables.
\end{abstract}


\section{Introduction.}

Let $\varphi : \mathbb C\longrightarrow \mathbb R$ be a $\mathcal C^2$-weight function
and consider the Hilbert spaces
$$L^2_{\varphi }=\{ f: \mathbb C \longrightarrow \mathbb C \  \mbox{measureable} \ : \
\|f\|^2_{\varphi }:=\int_{\mathbb C}|f(z)|^2 \,e^{-2\varphi (z)}\,d\lambda (z)<\infty \}.$$

It is essentially due to L. H\"ormander \cite{H}
that for a suitable weight function $\varphi $ and for every $f\in L^2_{\varphi }$
there exists $u\in L^2_{\varphi }$ satisfying $\ovprt u=f.$
In fact  there exists a continuous solution operator $\tilde S:L^2_{\varphi}
\longrightarrow L^2_{\varphi }$ for $\ovprt ,$ i.e. $\|\tilde S(f) \|_{\varphi }\le C
\|f\|_{\varphi }$ and $\ovprt \tilde S (f)=f,$ see also \cite{Ch} .

Let $A^2_{\varphi }$ denote the space of entire functions belonging to $L^2_{\varphi }$ and
let 
$$P_{\varphi }:L^2_{\varphi }\longrightarrow A^2_{\varphi }$$ 
denote the Bergman projection.
Then $S=(I-P_{\varphi })\tilde S $ is the uniquely determined canonical solution operator to $\ovprt ,$
i.e. $\ovprt S(f)=f $ and $S(f)\perp A^2_{\varphi }.$

In this paper we discuss the compactness of the canonical solution operator to $\ovprt $ on weighted
$L^2$-spaces. The question of compactness of the solution operator to $\ovprt $ is of interest for various
reasons - see \cite{FS1} and \cite{FS2} for an excellent survey and \cite{C}, \cite{CD},
\cite{K}, \cite{L}.

A similar situation appears in \cite{SSU} where the Toeplitz $C^*$ -algebra
$\mathcal T (\Omega )$ is considered and the relation between the structure of
$\mathcal T (\Omega )$ and the $\ovprt $-Neumann problem is discussed
(see \cite{SSU} ).

The connection of $\ovprt $ with the theory of Schr\"odinger operators with magnetic fields
appears in \cite{Ch}, \cite{B} and \cite{FS3}. 

For the case of one complex variable we use results of Iwatsuka
(\cite{I}) to discuss compactness of the canonical solution operator to $\ovprt$ (see also \cite{HeMo}).

Multiple difficulties arise in the case of several complex variables, mainly because the
geometric structures underlying the analysis become much more complicated.
We try to point out the different situation and the obstructions which appear in the case
of several complex variables.

\section{Schr\"odinger operators with magnetic fields in one complex variable.}

A nonnegative Borel measure $\nu $ defined on $\mathbb C$ is said to be doubling if
there exists a constant $C$ such that for all $z\in \mathbb C$ and $r\in \mathbb R^+,$
$$\nu (B(z,2r))\le C \nu (B(z,r)).$$
$\mathcal D$ denotes the set of all doubling measures $\nu $ for which there exists
a constant $\delta $ such that for all $z\in \mathbb C,$
$$\nu (B(z,1))\ge \delta.$$
Let $\varphi : \mathbb C\longrightarrow \mathbb R$ be a subharmonic function.
Then $\Delta \varphi $ defines a nonnegative Borel measure, which is finite on compact sets.

Let $\mathcal W$ denote the set of all subharmonic $\mathcal C^2$ functions $\varphi : \mathbb C
\longrightarrow \mathbb R$ such that $\Delta \varphi \in \mathcal D.$

\begin{thm}Let $\varphi \in \mathcal W.$ The canonical solution operator
$S:L^2_{\varphi}\longrightarrow L^2_{\varphi }$ to $\ovprt $ is compact if and only if
there exists a real valued continuous function $\mu $ on $\mathbb C$ such that
$\mu (z)\to \infty$ as $|z|\to \infty $ and
$$\int_{\mathbb C}\mathcal S \phi (z) \, \overline{\phi (z)}\, d\lambda (z)
\ge \int_{\mathbb C}\mu (z)\, |\phi (z)|^2\,d\lambda (z)$$
for all $\phi \in \mathcal C^{\infty}_0(\mathbb C),$ where
$$\mathcal S = -\frac{\partial^2}{\partial z \partial \overline z}-\frac{\partial
\varphi }{\partial \overline z}\,\frac{\partial}{\partial z}+ \frac{\partial \varphi }
{\partial z}\,\frac{\partial}{\partial \overline z}+\left |\frac{\partial \varphi }
{\partial z} \right |^2+\frac{\partial^2 \varphi}{\partial z \partial \overline z}.$$
\end{thm}

\vskip 0.3cm

\begin{proof}
Consider the equation $\ovprt u=f$ for $f\in L^2_{\varphi }.$ The canonical solution
operator to $ \ovprt $ gives a solution with minimal $L^2_{\varphi }$-norm.
We substitute $v=u\,e^{-\varphi }$ and $g=f\,e^{-\varphi }$
and the equation becomes
$$\overline D v=g \ , \ \mbox{where} \ \overline D = e^{-\varphi }\, \frac {\partial }
{\partial \overline z}\, e^{\varphi }.$$
$u$ is the minimal solution to the $\ovprt $-equation in $L^2_{\varphi }$ if and only if
$v$ is the solution to $\overline D v=g$ which is minimal in $L^2(\mathbb C ).$

The formal adjoint of $\overline D$ is $D=-e^{\varphi }\frac{\partial}{\partial z}e^{
-\varphi }.$ As in \cite{Ch} we define $\mbox{Dom}(\overline D)=\{f\in L^2(\mathbb C)
\ : \ \overline D f\in L^2(\mathbb C ) \}$ and likewise for $D.$ Then $\overline D $ and
$D$ are closed unbounded linear operators from $L^2(\mathbb C)$ to itself.
Further we define $\mbox{Dom}(\overline D D)=\{ u\in \mbox{Dom}(D) \ : \
Du\in \mbox{Dom}(\overline D )\}$ and we define $\overline D D $ as $\overline D \circ
D$ on this domain. Any function of the form $e^{\varphi }\,g,$ with $g\in \mathcal C^2_0$
belongs to $\mbox{Dom}(\overline D D)$ and hence $\mbox{Dom}(\overline D D)$ is dense
in $L^2(\mathbb C).$ Since $\overline D=\frac{\partial }{\partial \overline z}+
\frac{\partial \varphi }{\partial \overline z}$ and $D=-\frac{\partial }{\partial z}+
\frac{\partial \varphi }{\partial z}$ we see that
$$\mathcal S = \overline D D=-\frac{\partial^2}{\partial z \partial \overline z}-\frac{\partial
\varphi }{\partial \overline z}\,\frac{\partial}{\partial z}+ \frac{\partial \varphi }
{\partial z}\,\frac{\partial}{\partial \overline z}+\left |\frac{\partial \varphi }
{\partial z} \right |^2+\frac{\partial^2 \varphi}{\partial z \partial \overline z}$$
$$=-\frac{1}{4}\, ((d-iA)^2 - \Delta \varphi ),$$
where $A=A_1\,dx+A_2\,dy=-\varphi_y\,dx + \varphi_x\,dy.$ Hence $\mathcal S =\overline D D$ is a Schr\"odinger
operator with electric potential $\Delta \varphi $ and with magnetic field $B=dA,$ (\cite{CFKS}).

Now let $\|u\|^2=\int_{\mathbb C}|u(z)|^2\,d\lambda (z)$ for $u\in L^2(\mathbb C)$ and
$$(u,v)=\int_{\mathbb C}u(z)\overline{v(z)}\,d\lambda (z)$$
denote the inner product of
$L^2(\mathbb C).$

In \cite{Ch} the following results are proved :
If $u\in \mbox{Dom}(D)$ and $Du\in \mbox{Dom}(\overline D),$ then
$$\|Du\|^2= (\overline D(Du),u).$$
$\overline D D$ is a closed operator and
$$\|u\| \le C\|\overline D D u\|$$
for all $u\in \mbox{Dom}(\overline D D).$ Moreover, for any $f\in L^2(\mathbb C)$
there exists a unique $u\in \mbox{Dom}(\overline D D)$ satisfying $\overline D D u=f.$
Hence $\mathcal S^{-1}=(\overline D D )^{-1}$ is a bounded operator on $L^2(\mathbb C).$

\vskip 0.3 cm

Now we claim that
the canonical solution operator $S: L^2_{\varphi }\longrightarrow
L^2_{\varphi }$ to $\ovprt $ is compact if and only if $\mathcal S^{-1} : L^2(\mathbb C)
\longrightarrow L^2(\mathbb C )$ is compact.

For this we remark that $v$ is the minimal solution to $\ovprt v=g$ in $L^2_{\varphi }$ if and only if
$u=v\,e^{-\varphi }$  is the minimal solution to $\overline D u =g\,e^{-\varphi }$
in $L^2(\mathbb C).$ Hence the canonical solution operator $S$ to $\ovprt $ is compact
if and only if the canonical solution operator to $\overline D u =f $ is compact.
By the above properties of the operators $D$ and $\overline D$ we have
$$\|D\mathcal S^{-1}f\|^2=(\overline D D(\overline D D)^{-1}f, (\overline D D)^{-1}f)=
(f,(\overline D D)^{-1}f)\le \|\mathcal S^{-1}\|\|f\|^2, $$
hence
$$\|D\mathcal S^{-1}f\| \le \|\mathcal S^{-1}\|^{1/2} \, \|f\|$$
and $T=D\mathcal S^{-1}$ is a bounded operator on $L^2(\mathbb C)$ with
$\overline D Tf=f$ and $Tf\perp \mbox{ker}\overline D ,$ which means that
$T$ is the canonical solution operator to $\overline D u=f.$ Since $\mathcal S^{-1}$
is a selfadjoint operator (see for instance \cite{I}) it follows that
$$\mathcal S^{-1}=T^*T.$$
Since $T$ is compact if and only if $T^*T$ is compact (see \cite{W}), our claim is proved.

To prove the theorem we use Iwatsuka's result (\cite{I}) that the operator $\mathcal S$ has compact resolvent
if and only if the condition in Theorem 2.1 holds.
\end{proof}

\begin{thm} If $\varphi (z)=|z|^2,$ then the canonical solution operator
$S:L^2_{\varphi}\longrightarrow L^2_{\varphi }$ to $\ovprt $
fails to be compact.
\end{thm}

\begin{proof} In our case the magnetic field $B$ is the form $B=dA=B(x,y) dx\wedge dy=
\Delta \varphi dx\wedge dy.$
Hence for $\varphi (z)=|z|^2$ we have $\Delta \varphi (z)=4 $ for each $z \in \mathbb C.$
Let $Q_w $ be the ball centered at $w$ with radius $1.$ Then
$$\int_{Q_w}(|B(x,y)|^2+\Delta \varphi (z))\,d\lambda (z) $$
is a constant as $|w|\to \infty ,$ so the assertion follows from \cite{I} Theorem 5.2.
\end{proof}

\begin{thm} Let $\varphi \in \mathcal W$ and suppose that $\Delta \varphi (z)\to \infty$
as $|z|\to \infty .$ Then the
canonical solution operator $S:L^2_{\varphi}\longrightarrow L^2_{\varphi }$ to $\ovprt $ is
compact.
\end{thm}

\begin{proof} Since in our case $|B(x,y)|=\Delta \varphi (z)\to \infty,$ as $|z|\to \infty $
the conclusion follows from the proof of Theorem 2.1 and \cite{AHS}, \cite{D} or \cite{I}.
\end{proof}

{\bf Remark.} In \cite{Has2} it shown that for $\varphi (z) =|z|^2$
even the restriction of the canonical solution operator
$S$ to the Fock space $A^2_{\varphi }$ fails to be compact and that for $\varphi (z)=|z|^m \ ,
\ m>2$ the restriction of $S$ to $A^2_{\varphi }$ fails to be Hilbert Schmidt.

\section{Several complex variables.}

In \cite{Sch} it is shown that
the restriction of the canonical solution operator
to the Fock space $A^2_{\varphi }$ fails to be compact, where
$$\varphi (z)=|z_1|^m + \dots + |z_n|^m, $$
for $ m\geq 2$ and $n\geq 2.$ Hence the canonical solution operator
cannot be compact on the corresponding $L^2$-spaces.

Here we investigate the solution operator on $L^2$-spaces and try to generalize
the method from above for several complex variables.

Let $\varphi : \mathbb C^n \longrightarrow \mathbb R $ be a $\mathcal C^2$-weight function and consider the space
$$L^2(\mathbb C^n , \varphi )=\{ f:\mathbb C^n \longrightarrow \mathbb C \ : \ \int_{\mathbb C^n}
|f|^2\, e^{-2\varphi}\,d\lambda < \infty \}$$
and the space $L^2_{(0,1)}(\mathbb C^n, \varphi )$ of $(0,1)$-forms with coefficients in
$L^2(\mathbb C^n , \varphi ).$

For $v\in L^2(\mathbb C^n )$ let
$$\overline D v = \sum_{k=1}^n \left( \frac{\partial v}{\partial \ovli z_k}+
\frac{\partial \varphi}{\partial \ovli z_k}\, v \right) \, d\ovli z_k$$
and for $g=\sum_{j=1}^n g_j\, d\ovli z_j \in L^2_{(0,1)}(\mathbb C^n ) $ let
$$\ovli D^* g = \sum_{j=1}^n \left( \frac{\partial \varphi}{\partial z_j}\, g_j
-\frac{\partial g_j}{\partial z_j} \right) ,$$
where the derivatives are taken in the sense of distributions.
It is easy to see that $\ovprt u =f$ for $u\in L^2(\mathbb C^n , \varphi )$ and
$f\in L^2_{(0,1)}(\mathbb C^n, \varphi )$ if and only if $\ovli D v = g, $ where
$v= u\, e^{-\varphi }$ and $g= f\, e^{-\varphi }.$ It is also clear that
the necessary condition $\ovprt f=0$ for solvability holds if and only if
$\ovli D g =0 $ holds. Here
$$\ovli D g = \sum_{j,k=1}^n \left ( \frac{\partial g_j}{\partial \ovli z_k}
+\frac{\partial \varphi }{\partial \ovli z_k}\,g_j \right ) \, d\ovli z_k \wedge
d\ovli z_j.$$

Then

$$
\ovli D \, \ovli D^* g  =  \ovli D \left( \sum_{j=1}^n \left( \frac{\partial \varphi}{\partial z_j}\, g_j
-\frac{\partial g_j}{\partial z_j} \right) \right)$$
$$ =  \sum_{k=1}^n \left[ \sum_{j=1}^n \left(
\frac{\partial^2 \varphi}{\partial z_j \partial \ovli z_k}\, g_j
-\frac{\partial^2 g_j}{\partial z_j \partial \ovli z_k}+
\frac{\partial g_j}{\partial \ovli z_k}\, \frac{\partial \varphi}{\partial z_j}-
\frac{\partial g_j}{\partial z_j}\, \frac{\partial \varphi}{\partial \ovli z_k}+
\frac{\partial \varphi}{\partial z_j}\, \frac{\partial \varphi}{\partial \ovli z_k}\, g_j \right) \right] \,
d\ovli z_k .
$$

\begin{prp}
The operator $\ovli D \, \ovli D^*$ defined on $\mbox{Dom} \ovli D^* \cap \mbox{ker} \ovli D$ has the
form

$$  \sum_{k=1}^n \left[ \sum_{j=1}^n \left(
2\frac{\partial^2 \varphi}{\partial z_j \partial \ovli z_k}\, g_j
-\frac{\partial^2 \varphi}{\partial z_j \partial \ovli z_j}\, g_k
-\frac{\partial^2 g_k}{\partial z_j \partial \ovli z_j}
\right. \right.$$
$$\left. \left. +\frac{\partial g_k}{\partial \ovli z_j}\, \frac{\partial \varphi}{\partial z_j}-
\frac{\partial g_k}{\partial z_j}\, \frac{\partial \varphi}{\partial \ovli z_j}+
\frac{\partial \varphi}{\partial z_j}\, \frac{\partial \varphi}{\partial \ovli z_j}\, g_k \right) \right] \,
d\ovli z_k .
$$
\end{prp}

\begin{proof} The condition $\ovli D g=0$ means that
$$\frac{\partial g_j}{\partial \ovli z_k}+ \frac{\partial \varphi}{\partial \ovli z_k}\,g_j=
\frac{\partial g_k}{\partial \ovli z_j}+ \frac{\partial \varphi}{\partial \ovli z_j}\,g_k ,
$$
for $j,k=1,\dots , n.$ Now we apply the differentiation
$\frac{\partial }{\partial z_j}$ on both sides and obtain
$$\frac{\partial^2 \varphi}{\partial z_j \partial \ovli z_k}\, g_j
+\frac{\partial^2 g_j}{\partial z_j \partial \ovli z_k} +
\frac{\partial g_j}{\partial z_j}\, \frac{\partial \varphi}{\partial \ovli z_k}=
\frac{\partial^2 \varphi}{\partial z_j \partial \ovli z_j}\, g_k
+\frac{\partial^2 g_k}{\partial z_j \partial \ovli z_j} +
\frac{\partial g_k}{\partial z_j}\, \frac{\partial \varphi}{\partial \ovli z_j}$$

Using this for the formula for
$\ovli D \, \ovli D^* $ we get

$$
\ovli D \, \ovli D^* g =  \ovli D \left( \sum_{j=1}^n \left( \frac{\partial \varphi}{\partial z_j}\, g_j
-\frac{\partial g_j}{\partial z_j} \right) \right) $$

$$ =  \sum_{k=1}^n \left[ \sum_{j=1}^n \left(
2\frac{\partial^2 \varphi}{\partial z_j \partial \ovli z_k}\, g_j
-\frac{\partial^2 \varphi}{\partial z_j \partial \ovli z_j}\, g_k
-\frac{\partial^2 g_k}{\partial z_j \partial \ovli z_j}\right. \right. $$
$$\left. \left. +\frac{\partial g_k}{\partial \ovli z_j}\, \frac{\partial \varphi}{\partial z_j}-
\frac{\partial g_k}{\partial z_j}\, \frac{\partial \varphi}{\partial \ovli z_j}+
\frac{\partial \varphi}{\partial z_j}\, \frac{\partial \varphi}{\partial \ovli z_j}\, g_k \right) \right] \,
d\ovli z_k .
$$
\end{proof}

\begin{rem} The only term where $g_j$ appears in the last line is
$$2\frac{\partial^2 \varphi}{\partial z_j \partial \ovli z_k}\, g_j,$$
and we will get a diagonal system if we restrict to weight functions of a special form,
for instance
$\varphi (z)=|z_1|^2 + \dots +|z_n|^2 ,$
the case of the Fock space.
\end{rem}

\begin{prp}
Suppose that the weight function $\varphi $ is of the form
$$\varphi (z_1, \dots , z_n)=\varphi_1(z_1) + \dots + \varphi_n(z_n),$$
where $\varphi_j :\mathbb C \longrightarrow \mathbb R$ are $\mathcal C^2 $-functions
for $j=1,\dots , n.$

Then the equation $\ovli D \, \ovli D^* g = h ,$ for $h=\sum_{k=1}^n h_k \, d\ovli z_k,$
splits into the $n$-equations

$$2 \, \frac{\partial^2 \varphi}{\partial z_k \partial \ovli z_k}\, g_k +
\sum_{j=1}^n \left(
-\frac{\partial^2 \varphi}{\partial z_j \partial \ovli z_j}\, g_k
-\frac{\partial^2 g_k}{\partial z_j \partial \ovli z_j}\right.  $$
$$\left. +\frac{\partial g_k}{\partial \ovli z_j}\, \frac{\partial \varphi}{\partial z_j}-
\frac{\partial g_k}{\partial z_j}\, \frac{\partial \varphi}{\partial \ovli z_j}+
\frac{\partial \varphi}{\partial z_j}\, \frac{\partial \varphi}{\partial \ovli z_j}\, g_k \right) = h_k,$$

for $k=1, \dots , n.$
These equations can be represented as Schr\"odinger operators $\mathcal S_k $ with
magnetic fields, where

$$\mathcal S_k v =2 \, \frac{\partial^2 \varphi}{\partial z_k \partial \ovli z_k}\, v $$
$$+\sum_{j=1}^n \left(
-\frac{\partial^2 \varphi}{\partial z_j \partial \ovli z_j}\, v
-\frac{\partial^2 v}{\partial z_j \partial \ovli z_j}+
\frac{\partial v}{\partial \ovli z_j}\, \frac{\partial \varphi}{\partial z_j}-
\frac{\partial v}{\partial z_j}\, \frac{\partial \varphi}{\partial \ovli z_j}+
\frac{\partial \varphi}{\partial z_j}\, \frac{\partial \varphi}{\partial \ovli z_j}\, v
\right) $$

and $v$ is a $\mathcal C^2$-function.
The operators $\mathcal S_k$ can be written in the form
$$\mathcal S_k = \frac{1}{4} \left [
-\sum_{j=1}^n \left ( \frac{\partial }{\partial x_j}-ia_j \right )^2
-\sum_{j=1}^n \left ( \frac{\partial }{\partial y_j}-ib_j \right )^2 \right ]
+ V_k ,$$
where $z_j=x_j +iy_j$ and $a_j=-\frac{\partial \varphi }{\partial y_j} \ , \
b_j=\frac{\partial \varphi }{\partial x_j} ,$ for $j=1, \dots , n$ and
$$V_k=2\frac{\partial^2 \varphi}{\partial z_k \partial \ovli z_k}-
\sum_{j=1}^n \frac{\partial^2 \varphi}{\partial z_j \partial \ovli z_j},$$
for $k=1,\dots , n.$
\end{prp}

\begin{rem} 
If the weight function $\varphi $ is of the form
$$\varphi (z_1, \dots , z_n)=\varphi_1(z_1) + \dots + \varphi_n(z_n),$$
where $\varphi_j :\mathbb C \longrightarrow \mathbb R$ are $\mathcal C^2 $-functions
for $j=1,\dots , n,$ then the magnetic field of the Schr\"odinger operators
$\mathcal S_k$ is the 2-form
$$B = \sum_{j<l}B_{jl}\, d\tilde x_j \wedge d\tilde x_l,$$
where $\tilde x_{2j-1}=x_j$ , $\tilde x_{2j}=y_j$ , $\tilde a_{2j-1}=a_j$ ,
$\tilde a_{2j} =b_j$ for $j=1, \dots , n$ and
$$B_{jl}=\frac{1}{4}\, \left ( \frac{\partial \tilde a_l}{\partial \tilde x_j}
-\frac{\partial \tilde a_j}{\partial \tilde x_l} \right ).$$
If we write
$$|B|=\left ( \sum_{j<l} |B_{jl}|^2 \right )^{1/2},$$
\vskip 0.3 cm
then the assumptions on the weight function $\varphi $ imply that
\vskip 0.3 cm
$$|B|= \frac{1}{4} \, \left [ \sum_{j=1}^n \left ( \frac{\partial^2 \varphi}{\partial x_j^2}
+ \frac{\partial^2 \varphi}{\partial y_j^2} \right )^2 \right ]^{1/2} .$$
\vskip 0.5 cm
The electric potentials $V_k$ have the form
$$V_k=2\, \frac{\partial^2 \varphi}{\partial z_k \partial \ovli z_k}-
\sum_{j=1}^n \frac{\partial^2 \varphi}{\partial z_j \partial \ovli z_j}$$
$$=\frac{1}{2} \, \left ( \frac{\partial^2 \varphi}{\partial x_k^2}
+ \frac{\partial^2 \varphi}{\partial y_k^2}\right ) - \frac{1}{4}\,
\sum_{j=1}^n \left ( \frac{\partial^2 \varphi}{\partial x_j^2}+
\frac{\partial^2 \varphi}{\partial y_j^2} \right ).$$

Hence the the socalled effective potentials (see \cite{KS} , Corollary 1.14)
$$V_{k, \mbox{eff}}^{\delta}= V_k + \frac{\delta}{n-1}\, |B| \ , \quad \delta \in [0,1)$$
do not tend to infinity as $|z|$ tends to infinity for weight functions like
$$\varphi (z)= \sum_{j=1}^n |z_j|^2 ,$$
causing the obstructions for the Schr\"odinger operators $\mathcal S_k$ to have compact resolvents.

\end{rem}

{\scshape
\begin{flushright}
\begin{tabular}{l}
Fakult\"at f\"ur Mathematik, Universit\"at Wien,\\
Nordbergstraße 15, \\
A-1090 Wien, \\
Austria \\
{\upshape e-mail: friedrich.haslinger@univie.ac.at}\\
\end{tabular}
\end{flushright}
}

\end{document}